\documentclass[11pt]{article}

\usepackage{amssymb, amsmath, amsthm, tikz, xspace,subfigure, graphicx}
\usepackage[left=1in,top=1in,right=1in]{geometry}
\usetikzlibrary{positioning, shapes.geometric, calc, intersections}
\tikzstyle{vertex}=[circle,fill=black,inner sep=2pt]
\tikzstyle{vertrect}=[draw,rectangle,inner sep=2pt]
\tikzstyle{vertdia}=[draw,diamond,inner sep=2pt]

\date{}

\theoremstyle{plain}
      \newtheorem{theorem}{Theorem}[section]

      \newtheorem{problem}[theorem]{Problem}
      
      \newtheorem{corollary}[theorem]{Corollary}
      
      \newtheorem{conjecture}[theorem]{Conjecture}
\theoremstyle{definition}
      \newtheorem{definition}[theorem]{Definition}

\theoremstyle{remark}

	\newcommand{\ZZ}{{\mathbb Z}}

\def\twr{\mbox{\rm twr}}

\title{A survey of  hypergraph Ramsey problems}

\author{Dhruv Mubayi\thanks{Department of Mathematics, Statistics, and Computer Science, University of Illinois, Chicago, IL, 60607 USA.  Research partially supported by NSF grant DMS-1300138. Email: {\tt mubayi@uic.edu}} \and Andrew Suk\thanks{Department of Mathematics,  University of California at San Diego, La Jolla, CA, 92093 USA. Supported by NSF grant DMS-1800736, an NSF CAREER award, and an Alfred Sloan Fellowship. Email: {\tt asuk@ucsd.edu}.}}

\begin{document}

\maketitle

\begin{abstract}
 The classical \emph{hypergraph Ramsey number} $r_k(s,n)$ is the minimum $N$ such that for every red-blue coloring of the $k$-tuples of $\{1,\ldots, N\}$, there are $s$ integers such that every $k$-tuple among them is red, or $n$ integers such that every $k$-tuple among them is blue. We survey a variety of problems and results in hypergraph Ramsey theory that have grown out of understanding the quantitative aspects of  $r_k(s,n)$. Our focus is on recent developments and open problems.
\end{abstract}

\section{Introduction}

 A $k$-uniform hypergraph $H$ ($k$-graph for short) with vertex set $V$ is a collection of $k$-element subsets of $V$.  We write $K^{(k)}_n$ for the complete $k$-graph on an $n$-element vertex set.  The \emph{Ramsey number} $r_k(s,n)$ is the minimum $N$ such that every red-blue coloring of the edges of  $K^{(k)}_N$  contains a monochromatic red copy of $K_s^{(k)}$ or a monochromatic blue copy of $K^{(k)}_n$.  The existence of $r_k(s,n)$ follows from the celebrated theorem of Frank Ramsey from 1930 \cite{Ram}.  However, the asymptotic behavior of $r_k(s,n)$ is still not well understood.

  In this survey we focus on open problems and results related to generalizations and extensions of $r_k(s,n)$ in the {\em hypergraph} case, i.e., when $k \ge 3$
  (we refer the reader to \cite{CFSsurvey} for a survey of Graph Ramsey theory). Our emphasis is on recent results and although we believe we have touched on most important developments in this area, this survey is not an exhaustive compendium of all work in hypergraph Ramsey theory.

\section{General notation}

The full statement of Ramsey's theorem extends to multiple colors and to general hypergraphs as follows.   Given an integer $q\geq 2$ and $k$-uniform hyergraphs $H_1,\ldots, H_q$, there is a minimum $r_k(H_1,\ldots, H_q) = N$, such that every $q$-coloring of the edges of $K^{(k)}_N$ contains a copy of $H_i$ in the $i$th color.  In the special case that $H = H_1 = \cdots = H_q$, we simply write

$$r_k(H; q) = r_k(\underbrace{H,\ldots, H}_{q\textnormal{ times}}).$$  If $H_i = K^{(k)}_{n_i}$, we use the simpler notation $r_k(n_1,\ldots, n_q)$ and $r_k(n;q) = r_k(\underbrace{n,\ldots, n}_{q\textnormal{ times}}).$

\section{Diagonal Ramsey numbers}

 \emph{Diagonal Ramsey numbers} refer to the special case when $s = n$, i.e. $r_k(n,n)$, and have been studied extensively over the past 80 years.  Classic results of Erd\H os and Szekeres \cite{ES35} and Erd\H os \cite{E47} imply that $2^{n/2} < r_2(n,n) \leq 2^{2n}$ for every integer $n > 2$.  While small improvements have been made in both the upper and lower bounds for $r_2(n,n)$ (see \cite{S75,C09}), the constant factors in the exponents have not changed over the last 70 years.

Unfortunately for 3-graphs, our understanding of $r_3(n,n)$ is much less than in the graph case.  A result of Erd\H os, Hajnal, and Rado \cite{EHR} gives the best known lower and upper bounds for $r_3(n,n)$, $$2^{c_1n^2}<r_3(n,n)<2^{2^{c_2n}},$$ where $c_1$ and $c_2$ are absolute constants.  Another proof of the lower bound above was given by Conlon, Fox, and Sudakov in \cite{CFS}, which will be discussed in more detail in Section \ref{offdiagsec}.  For $k \geq 4$, there is also a difference of one exponential between the known lower and upper bounds for $r_k(n,n)$, that is, \begin{equation}\label{diag}\twr_{k-1}(c_1n^2) \leq r_k(n,n) \leq \twr_k(c_2n),\end{equation}

\noindent where the \emph{tower function} $\twr_k(x)$ is defined by $\twr_1(x) = x$ and $\twr_{i + 1}(x) = 2^{\twr_i(x)}$ (see \cite{ES35,ER,EH72}).  A notoriously difficult conjecture of Erd\H os, Hajnal, and Rado states that the upper bound in (\ref{diag}) is essentially the truth, that is, there are constructions which demonstrate that $r_k(n,n) > \twr_k(cn)$, where $c = c(k)$.  The crucial case is when $k = 3$, since a double exponential lower bound for $r_3(n,n)$ would verify the conjecture for all $k\geq 4$ by using the well-known stepping-up lemma of Erd\H os and Hajnal (see \cite{GRS}).

\begin{conjecture}[Erd\H os]\label{3conj}
For $n\geq 4$ we have $r_3(n,n) > 2^{2^{cn}}$, where $c$ is an absolute constant.

\end{conjecture}

It is worth mentioning that Erd\H os offered a \$500 reward for a proof of this conjecture (see~\cite{Chung}), and his conjecture is supported by the fact that a double exponential lower bound is known if one allows four colors. More precisely, Erd\H os and Hajnal (see \cite{GRS}) showed that $r_3(n;4) > 2^{2^{cn}}$, and for three colors, the best known lower bound for $r_3(n; 3)$ is due to Conlon, Fox, and Sudakov \cite{CFS} who showed that $r_3(n; 3) > 2^{n^{c\log n}}$.  There is some evidence that perhaps Conjecture \ref{3conj} is false, and we refer the interested reader to \cite{almostmono,CFR} for two results in this direction.

\section{Off-Diagonal Ramsey numbers}\label{offdiagsec}

 \emph{Off-diagonal Ramsey numbers}, $r_k(s,n)$, refer to the special case when $k,s$ are fixed and $n$ tends to infinity.  It is known \cite{AKS,Kim,B,BK} that $r_2(3,n) =\Theta(n^2/\log n)$, and more generally for fixed $s > 3$, $r_2(s,n) = n^{\Theta(1)}$.  For 3-graphs, Conlon, Fox and Sudakov~\cite{CFS} proved that there are absolute constants $c, c'>0$ such that for all $4 \le s \le n$, 
$$2^{csn\log \left(\frac{n}{s} + 1\right)} < r_3(s,n) <2^{ (c'n/s)^{s-2} \log (n/s)}.$$
For $s = n$, this gives another proof that $r_3(n,n) > 2^{cn^2}$.

   For $k$-graphs, where $s > k\geq 4$, it is known that
$r_k(s,n) \leq \twr_{k-1}(n^{c}),$ where $c = c(s)$ \cite{ER}.   Erd\H os and Hajnal proved that
 \begin{equation}\label{off}r_k(s,n)> \twr_{k-1}(c'n),\end{equation} for $k\geq 4$ and $s \geq 2^{k-1} - k + 3$, where $c' = c'(s)$.  They conjectured that a similar bound should hold  for smaller $s$ as follows.

\begin{conjecture}[Erd\H os-Hajnal~\cite{EH72}] \label{ehoffconj}
Fix $4\le k<s$. There are constants $c$ and $c'$ such that
$$\twr_{k-1}(cn) < r_k(s,n)< \twr_{k-1}(c'n).$$
\end{conjecture}
 Actually, this was part of a more general conjecture that they posed in that paper which will be discussed in Section 4.  Erd\H os and Hajnal (see \cite{GRS}) showed that $r_4(7,n) > 2^{2^{cn}}$, and the authors \cite{MS15} and Conlon, Fox, and Sudakov \cite{CFS15} independently verified the conjecture for $k\geq 4$ and $s \geq k + 3$ (using different constructions).  However, showing that $r_4(5,n)$ and $r_4(6,n)$ grows double exponentially in a power of $n$ seemed to be much more difficult.

Just as for diagonal Ramsey numbers, a double exponential in $n^c$ lower bound for $r_4(5,n)$ and $r_4(6,n)$ would imply $r_k(k + 1,n) > \twr_{k-1}(n^{c'})$ and $r_k(k + 2,n) > \twr_{k-1}(n^{c'})$ respectively, for all fixed $k\geq 5$, by a variant of the Erd\H os-Hajnal stepping up lemma. In \cite{MSr5r6}, the authors established the following lower bounds for $r_4(5,n)$ and $r_4(6,n)$, which represents the current best bounds:
for all $n\geq 6$,
$$r_4(5,n)> 2^{n^{c\log n}} \hspace{1cm}\textnormal{and}\hspace{1cm}r_4(6,n)> 2^{2^{cn^{1/5}}},$$
where $c > 0$ is an absolute constant. More generally, for $n > k \geq 5$, there is a $c = c(k) > 0$ such that
$$r_k(k+1,n)> \twr_{k-2}(n^{c\log n}) \hspace{1cm}\textnormal{and}\hspace{1cm}r_k(k+2,n)> \twr_{k-1}(cn^{1/5}).$$
A standard argument in Ramsey theory together with results in~\cite{CFS} for 3-graphs  yields  $$r_k(k+2, n) <
\twr_{k-1}(c'n^3\log n),$$ so we now know the tower growth rate of $r_k(k+2, n)$.  It remains an open problem to prove that $r_4(5,n)$ is double exponential in a power of $n$.

\begin{conjecture}\label{offconj}
For $n\geq 5$, there is an absolute constant $c> 0$ such that $r_4(5,n)  > 2^{2^{n^c}}$.
\end{conjecture}

In \cite{MS15}, the authors established a connection between diagonal and off-diagonal Ramsey numbers, by showing that a solution to Conjecture \ref{3conj} implies a solution to Conjecture~\ref{offconj} (see Section~\ref{ordsect} for more details).

\section{The Erd\H os-Hajnal Problem}
As mentioned in previous sections, it is a major open problem to determine if $r_3(n,n)$ and $r_4(5,n)$ grow double exponentially in a power of $n$.  In order to shed more light on these questions,  Erd\H os and Hajnal~\cite{EH72} in 1972 considered the  following more general parameter.

\begin{definition} For integers $2\le k < s <n$ and $2 \le t \le {s \choose k}$, let $r_k(s,t;n)$ be the minimum $N$ such that every red/blue coloring of  the edges of $K^{(k)}_N$ results in a monochromatic blue copy of $K_n^{(k)}$ or has a set of $s$ vertices which contains at least $t$ red edges.
\end{definition}
The function $r_k(s, t; n)$ encompasses several fundamental problems which have been studied for a while.  Clearly $r_k(s,n) = r_k(s, {s \choose k}; n)$ so $r_k(s,t;n)$ includes classical Ramsey numbers.  In addition to off-diagonal and diagonal Ramsey numbers already mentioned, the function $r_k(k +1, k + 1; k + 1)$ has been studied in the context of the Erd\H os-Szekeres theorem and Ramsey numbers of ordered tight-paths by several researchers~\cite{DLR, EM, FPSS, MS, MSW}, the more general function $r_k(k+1, k+1; n)$ is related to high dimensional tournaments~\cite{LM}, and even the very special case $r_3(4,3;n)$ has tight connections to quasirandom hypergraph constructions~\cite{BR, KNRS, LM1, LM2}.


 The main conjecture of Erd\H os and Hajnal \cite{EH72} for $r_k(s,t;n)$ is that, as $t$ grows from $1$ to ${s\choose k}$, there is a well-defined value $t_1=h_1^{(k)}(s)$ at which $r_k(s,t_1-1;n)$ is polynomial in $n$ while $r_k(s,t_1;n)$ is exponential in a power of $n$, another well-defined value $t_2=h_2^{(k)}(s)$ at which it changes from exponential to double exponential in a power of $n$ and so on,  and finally a well-defined value $t_{k-2}=h_{k-2}^{(k)}(s)<{s \choose k}$
  at which it changes from $\twr_{k-2}$ to $\twr_{k-1}$ in a power of $n$. They were not able to offer a conjecture as to what $h_i^{(k)}(s)$ is in general, except when $i=1$ and when $s=k+1$.

  $\bullet$ When $i=1$, they conjectured that $t_1=h_1^{(k)}(s)$ is one more than the number of edges in the $k$-graph obtained by taking a complete $k$-partite $k$-graph on $s$ vertices with almost equal part sizes, and repeating this construction recursively within each part. Erd\H os offered \$500 for a proof of this (see \cite{Chung}).

 $\bullet$ When $s=k+1$, they conjectured that $h_i^{(k)}(k+1)=i+2$, that is, $r_k(k + 1, 2;n)$ is polynomial in $n$, $r_k(k +1,3;n)$ is exponential in a power of $n$,  $r_k(k +1,4;n)$ is double exponential in a power of $n$, and etc.~such that at the end, both $r_k(k+1,k;n)$ and $r_k(k+1,k+1;n)$ are $\twr_{k-1}$ in a power of $n$.  They proved this for $i=1$ via the following:

\begin{theorem} [Erdos-Hajnal~\cite{EH72}] \label{eh}
For $k \ge 3$, there are  positive $c= c(k)$ and $c'=c'(k)$ such that
$$r_k(k+1, 2; n) < cn^{k-1} \qquad \hbox{ and } \qquad
r_k(k+1, 3; n) > 2^{c'n}.$$
\end{theorem}

Results of  R\"odl-\v Sinajov\'a \cite{RS} on partial Steiner systems, and of Kostochka-Mubayi-Verstra\"ete~\cite{KMV}
on independent sets in hypergraphs, determine the order of magnitude of the function $r_k(k+1,2;n)$ as follows.  For each $k \ge 3$ there exist positive $c=c_k$ and $c'=c'(k)$
 such that
 $$ c' n^{k-1}/\log n <  r_k(k+1,2;n) < c \, n^{k-1}/\log n.$$

  For the $t = 3$ case, the authors in \cite{MS16} showed that for $k\geq 3$, there are positive  $c = c(k)$ and $c' = c'(k)$ such that

\begin{equation}\label{ehk12}
2^{cn^{k-2}} \leq r_k(k + 1,3;n) \leq 2^{c'n^{k-2}\log n }.
\end{equation}

For general $t$, the methods of Erd\H os and Rado \cite{ER} show that there exists $c=c(k,t)>0$ such that \begin{equation}\label{ehre}r_k(k+1, t; n) \leq \twr_{t-1}(n^c),\end{equation} for $3 \leq t\leq k$.  Erd\H os and Hajnal conjectured that this upper bound is the correct tower growth rate for $ r_k(k+1, t; n)$.

\begin{conjecture} [Erdos-Hajnal~\cite{EH72}] \label{ehconj}
For $k\geq 3$ and $2 \le t \le k$, there exists  $c=c(k,t)>0$ such that
$$r_k(k+1, t; n) \geq \twr_{t-1}(c \, n).$$
\end{conjecture}

Note that when $t = k+1$, the results from the previous section states that $r_k(k+1,k+1;n) = r_k(k+1,n) \leq \twr_{k-1}(n^{c'})$ where $c' = c'(k,t)$.

Hence for 3-graphs, $r_3(4,t;n)$ is fairly well understood.  We know that $r_3(4,2;n)$ is polynomial in $n$, and both $r_3(4,3;n)$ and $r_3(4,4;n)$ are exponential in $n^{1 + o(1)}$.  See \cite{CFS} for more results on $h_1^{(3)}(s)$ for $s > 4$.  Unfortunately for 4-graphs, we do not have a good understanding of $r_4(5,t;n)$ when $4 \leq t\leq 5$.  The best known upper and lower bounds for $r_4(5,4;n)$ are obtained by (\ref{ehk12}) and (\ref{ehre}), which give $2^{cn^2} < r_4(5,4;n) < 2^{2^{n^{c}}}.$  Notice that Conjecture~\ref{ehconj} states that $r_4(5,4;n)$ grows double exponential in a power of $n$, but we don't even know if $r_4(5,5;n) = r_4(5,n)$ is double exponential in a power of $n$.  Likewise, for 5-graphs, not much is known about $r_5(6,4;n)$ and $r_5(6,5;n)$.  Combining (\ref{ehk12}) and (\ref{ehre}) gives

\begin{equation}\label{r52}
2^{c'n^3} < r_5(6,4;n) < 2^{2^{n^{c}}} \hspace{1cm}\textnormal{and}\hspace{1cm} 2^{c'n^3} < r_5(6,5;n) < 2^{2^{2^{n^{c}}}}.
\end{equation}

\begin{problem}

Determine the tower growth rate of $r_4(5,4;n)$, $r_5(6,4;n)$, and $r_5(6,5;n)$.

\end{problem}

    However for $k$-graphs, when $k\geq 6$, the authors in \cite{MS16} settled Conjecture \ref{ehconj} in almost all cases in a strong form, by determining the correct tower growth rate, and in half of the cases also determining the correct power of $n$ within the tower.

\begin{theorem}[Mubayi-Suk \cite{MS16}]\label{MSEH}
For $k\geq 6$ and $4 \le t \le k-2$, there are positive $c =c(k,t)$ and $c'=c'(k,t)$ such that
$$\twr_{t-1}(c' n^{k-t + 1}\log n) \, \ge \, r_k(k+1,t; \, n)
 \, \ge \,  \begin{cases}
\twr_{t-1}(c \, n^{k-t + 1}) \qquad  \hbox{ if $k-t$ is even}\\
\twr_{t-1}(c \, n^{(k-t + 1)/2}) \hskip8pt \hbox{ if $k-t$ is odd.}
\end{cases}
$$
\end{theorem}

When $k\geq 6$ and $t \in \{k-1, k\}$, Conjecture \ref{ehconj} remains open.  We note that the upper bound in Theorem~\ref{MSEH} also holds when $k-1 \leq t \leq k$.  The best known upper and lower bounds for $r_k(k+1,k-1;n)$ and $r_k(k+1,k;n)$, also due to the authors \cite{MS16}, are

$$\twr_{k-3}(c \, n^3) \le  r_k(k+1, k-1;\, n)  \le \twr_{k-2}(c' \, n^2),$$
and
$$\twr_{k-3}(c \, n^3) \le r_k(k+1, k;\, n)  \le \twr_{k-1}(c' \, n).$$

In fact, by using the stepping-up lemma established in \cite{MS16}, any improvement in the lower bound for $r_5(6,4;n)$ and $r_5(6,5;n)$ in (\ref{r52}) would imply a better lower bound for $ r_k(k+1, k-1;\, n)$ and $r_k(k+1, k;\, n)$ respectively.

\medskip

\section{The Erd\H os-Rogers Problem}
An $s$-independent set in a $k$-graph $H$ is a vertex subset that contains no copy of $K_s^{(k)}$. So if $s=k$, then it is just an independent set.  Let $\alpha_s(H)$ denote the size of the largest $s$-independent set in $H$.

\begin{definition}
 For $k \le s < t < N$, the Erd\H os-Rogers function $f^k_{s,t}(N)$ is the minimum of $\alpha_s(H)$ taken over all $K_t^{(k)}$-free $k$-graphs $H$ of order
 $N$.
\end{definition}

 To prove the lower bound $f_{s,t}^{k}(N)\ge n$, one must show that every $K_{t}^{(k)}$-free $k$-graph on $N$ vertices contains an $s$-independent set with
 $n$ vertices. On the other hand, to prove the upper bound $f_{s,t}^{(k)}(N) < n$, one must construct a $K_{t}^{(k)}$-free $k$-graph $H$ of order $N$ with
 $\alpha_s(H) < n$.

The problem of determining $f_{s,t}^{k}(n)$ extends that of finding Ramsey numbers. Formally,
$$
r_k(s,n) = \min \{ N : f_{k,s}^{k}(N) \ge n\}.
$$

For $k=2$, the above function was first considered by Erd{\H o}s and Rogers~\cite{ERog} only for $t=s+1$, which is perhaps the most interesting case. So in this case we wish to construct a $K_{s+1}$-free graph on $N$ vertices such that the $s$-independence number is as small as possible.
Since then the function has been studied by several researchers culminating in the work of Wolfowitz~\cite{Wo} and Dudek, Retter and R\"odl~\cite{DRR} who proved the upper bound that follows (the lower bound is due to Dudek and the first author~\cite{DM}): for every $s\ge 3$ there are positive constants $c_1$ and $c_2=c_2(s)$ such that
\begin{equation} \label{fs}
c_1\left(\frac{ N \log N }{\log\log N}\right)^{1/2}< f^2_{s,s+1}(N)< c_2 (\log N)^{4s^2}N^{1/2}.\end{equation}
 The problem of estimating the Erd\H os-Rogers function for $k>2$ appears to be much harder. Let us denote $$g(k,N)=f^k_{k+1, k+2}(N).$$
In other words, $g(k,N)$ is the minimum $n$ such that every $K_{k+2}^{(k)}$-free $k$-graph on $N$ vertices has the property that every $n$-set of vertices has a copy of $K_{k+1}^{(k)}$.
 With this notation, the bounds in (\ref{fs}) for $s=3$ imply that $g(2,N)=N^{1/2+o(1)}$.

 Dudek and the first author \cite{DM} proved that $(\log N)^{1/4+o(1)} < g(3,N) < O(\log N)$, and more generally, that there are positive  $c_1=c_1(k)$ and $c_2=c_2(k)$ with
\begin{equation} \label{1} c_1( \log_{(k-2)}N)^{1/4} < g(k,N) < c_2(\log N)^{1/(k-2)},\end{equation}
where $\log_{(i)}$ is the log function iterated $i$ times.
The exponent 1/4 in (\ref{1}) was improved to 1/3 by Conlon, Fox and Sudakov \cite{CFS14}.
Both sets of authors asked whether the upper bound could be improved (presumably  to an iterated log function). This was achieved by the current authors~\cite{MS15c} who proved that  
for $k \ge 14$,
 $$g(k,N) = O( \log_{(k-13)} N ).$$
 It remains an open problem to determine the correct number of iterations (which may well be $k-2$).
 We pose this as a conjecture.
 \begin{conjecture}
For all $k \ge 3$, there are $c_1, c_2>0$ such that
$$c_1 \log_{(k-2)} N < g(k,N) < c_2 \log_{(k-2)} N .$$
\end{conjecture}

\section{The Erd\H os-Gy\'arf\'as-Shelah Problem}
A $(p,q)$-coloring of $K_N^{(k)}$ is an edge-coloring of $K_N^{(k)}$ that gives every copy of $K_p^{(k)}$ at least $q$ colors. Let $f_k(N,p,q)$ be the minimum number of colors in a $(p,q)$-coloring of $K_N^k$.

The problem of determining $f_k(N,p,q)$ for fixed $k,p,q$ has a long history, beginning with its introduction by Erd\H os and Shelah~\cite{E1, E2}, and  subsequent investigation (for graphs) by Erd\H os and Gy\'arf\'as~\cite{EG}.  Since
$$f_k(N,p,2)=t \qquad  \Longleftrightarrow \qquad r_k(p; t)\ge N+1 \quad \hbox{ and }\quad  r_k(p; t-1)\le N,$$
most of the effort on determining $f_k(N, p,q)$ has been for $q>2$.
As mentioned above, Erd\H os and Gy\'arf\'as~\cite{EG} initiated a systematic study of this parameter for graphs and posed many open problems. One main question was to determine the minimum $q$ such that
$f_2(N, p,q)= N^{o(1)}$ and $f_2(N, p,q+1)> N^{c_p}$ for some $c_p>0$. For $p=3$ this value is clearly $q=2$ as $f_2(N, 3, 2)= O(\log N)$ due to the easy bound $r_2(3; t)>2^t$,
 while $f_2(N, 3, 3) =\chi'(K_N) \ge N-1$. Erd\H os and Gy\'arf\'as proved that $f_2(N, p, p)>N^{c_p}$
 and asked whether  $f_2(N, p, p-1)= N^{o(1)}$. The first open case was $f_2(N,4,3)$, which was shown to be $N^{o(1)}$ by the first author~\cite{M} and later $\Omega(\log N)$ (see~\cite{FS, KM}).  The same upper bound was shown for $f(N,5,4)$ in~\cite{EMub}. Conlon, Fox, Lee and Sudakov~\cite{CFLS2} recently extended this construction considerably by proving  that
$f_2(N,p,p-1)=N^{o(1)}$ for all fixed $p \ge 4$. Their result is sharp in the sense that $f_2(N,p,p)=\Omega(N^{1/(p-2)})$. The exponent $1/(p-2)$ was shown to be sharp for $p=4$ by the first author~\cite{M2} and recently also for $p=5$ by Cameron and Heath~\cite{CH} via explicit constructions.

The first nontrivial hypergraph case is $f_3(N,4,3)$ and this function has tight connections to Shelah's breakthrough proof~\cite{S} of primitive recursive bounds for the Hales-Jewett numbers. Answering a question of Graham, Rothschild and Spencer~\cite{GRS}, Conlon, Fox, Lee and Sudakov showed that
  $$f_3(N,4,3)= N^{o(1)}.$$
They also posed a variety of basic questions about $f_k(N,p,q)$ and related parameters including the following generalization of the Erd\H
 os-Gy\'arf\'as problem for hypergraphs. Using a variant of the pigeonhole argument for hypergraph Ramsey numbers due to Erd\H os and Rado,~\cite{CFLS} proved that
 $$f_k\left(N, p, {p-i \choose k-i} +1\right) = \Omega(\log_{(i-1)}N^{c_{p,k,i}})$$
 where $\log_{(0)}(x)=x$ and, as usual, $\log_{(i+1)}x=\log \log_{(i)}x$ for $i \ge 0$.

 \begin{problem} [Conlon-Fox-Lee-Sudakov~\cite{CFLS}] \label{cflsprob}
For $p>k \ge 3$ and $0<i<k$ prove that $f_k(N, p, {p-i \choose k-i})$ is substantially smaller than $f_k(N, p, {p-i \choose k-i}+1)$, in particular, prove that $f_k(N, p, {p-i \choose k-i})$ is much smaller than $\log_{(i-1)}N$.
\end{problem}

One natural way to interpret  this problem is that it asks  whether
$$f_k\left(N, p, {p-i \choose k-i}\right)=(\log_{(i-1)}N)^{o(1)}?$$ The case $k=2$ is precisely the Erd\H os-Gy\'arf\'as problem and the case
$k=3, p=4, i=1$ is to prove that $f_3(N,4,3)= N^{o(1)}$ which was established in~\cite{CFLS}.
The next open case is $k=3, p=5, i=2$, which asks whether $f_3(N,5,3) = (\log N)^{o(1)}$.
This was solved with a better bound by the first author~\cite{Mlocal}, who showed that $$f_3(N,5,3)= e^{O(\sqrt{\log\log N})}= (\log N)^{O(1/\sqrt{\log\log N})}.$$  No other nontrivial cases of Problem~\ref{cflsprob} have been solved. We refer the reader to~\cite{CFLS} for related problems and results.

\section{More off-diagonal problems}

In this section we consider $k$-graph Ramsey numbers of the form $r_k(H, n):= r_k(H, K_n^{(k)})$
where $H$ is a (fixed) $k$-graph and $n$ grows.

 \subsection{$K_4^{(3)}$ minus an edge and a generalization}

Let $K_4^{(3)}\setminus e$ denote the 3-graph on four vertices, obtained by removing one edge from $K_4^{(3)}$.  A simple argument of Erd\H os and Hajnal \cite{EH72}  implies $r(K_4^{(3)}\setminus e,K_n^{(3)}) < (n!)^2$.   This was generalized in~\cite{MS15c} as follows.  A \emph{$k$-half-graph}, denote by $B=B^{(k)}$, is a $k$-graph on $2k-2$ vertices, whose vertex set is of the form $S\cup T$, where $|S| = |T| = k-1$, and whose edges are all $k$-subsets that contain $S$, and one $k$-subset that contains $T$.  So $B^{(3)}=K_4^{(3)}\setminus e$. Write $r_k(B, n)=r(B^{(k)}, K_n^{(k)})$.  It was shown in \cite{MS15c} that for each $k\geq 4$ there exists $c=c_k$ such that
$$2^{cn} < r_k(B,n) < (n!)^{k-1}.$$

A problem that goes back to the 1972 paper of Erd\H os and Hajnal (for $k=3$) is to improve the lower bound above. Indeed, $r_3(B, n)=r_3(4, 3; n)$ and this is therefore a very special case of the Erd\H os-Hajnal problem discussed earlier.

\begin{problem}
Show that for each $k\geq 3$ there exists $c=c_k$ such that
$r_k(B, n) > 2^{c n \log n}$.
\end{problem}

\subsection{Independent neighborhoods}

\begin{definition} A $k$-uniform triangle $T^{(k)}$ is a set of
$k+1$ edges $b_1, \ldots, b_k, a$ with $b_i \cap b_j=R$ for all $i<j$ where $|R|=k-1$ and $a = \cup_i (b_i-R)$. In other words, $k$ of the edges share a common $(k-1)$-set of vertices, and the last edge contains the remaining point in all these previous edges.
\end{definition}

 When $k=2$, then $T^{(2)}=K_3$, so in this sense $T^{(k)}$ is a generalization of a graph triangle.  We may view a $T^{(k)}$-free $k$-graph
as one in which all  neighborhoods are independent sets, where the neighborhood of an $R \in {V(H)\choose k-1}$ is
$\{x: R \cup \{x\} \in H\}$.  As usual, write $r_k(T, n)$ for $r(T^{(k)}, K_n^{(k)})$.

Bohman, Frieze and Mubayi~\cite{BFM} proved that for fixed $k \ge 2$, there are positive constants $c_1$ and $c_2$ with
$$c_1\frac{n^k}{(\log n)^{k/(k-1)}}<r_k(T, n)< c_2 n^k.$$
They conjectured that the upper bound could be improved to $o(n^k)$ and believed that the log factor in the lower bound could also be improved. Results of
 Kostochka-Mubayi-Verstra\"ete~\cite{KMV} proved this  and then Bohman-Mubayi-Picollelli~\cite{BMP} achieved a matching lower bound by analyzing the hypergraph independent neighborhood process. This may be viewed as a hypergraph generalization of the results of Ajtai-Koml\'os-Szemer\'edi~\cite{AKS} for graphs.
\begin{theorem} [Kostochka-Mubayi-Verstra\"ete~\cite{KMV}, Bohman-Mubayi-Picollelli~\cite{BMP}]
For fixed $k \ge 3$ there are positive constants $c_1$ and $c_2$ with
$$c_1\frac{n^k}{\log n}<r_k(T, n)< c_2 \frac{n^k}{\log n}.$$
\end{theorem}

\subsubsection{Unordered tight-paths versus cliques}
An (unordered) 3-uniform tight-path $TP^{(3)}_s = TP_s$ is the $3$-graph with vertex set $\{v_1,\ldots, v_s\}$ and
edge set $\{\{v_i, v_{i+1}, v_{i+2}\}: i \in \{1,\ldots, s-2\}\}$.  Note that the vertex set $\{v_1,\ldots, v_s\}$ is not ordered.
Results of Phelps and R\"odl~\cite{phelpsrodl} imply that there are
 $c_1$ and $c_2$ such that
\[
c_1 n^2/\log n < r_3(TP_4, n) < c_2 n^2/\log n.
\]
It is easy to prove that for all $s \ge 5$, there is $c=c_s$ such that
$r_3(TP_s, n) < c\, n^2$. A matching lower bound for $s\geq 6$ was provided by Cooper and Mubayi \cite{CMsparse} with the following construction.  Let $H$ be a 3-graph where $V(H) = [n]\times [n]$, and $E(H)  = \{\{ab,ac,db\} \in [n]\times [n]: c > b, d > a\}$.  It is easy to see that $H$ is $TP_6$-free and $\alpha(H) < 2n$.  Thus, For $s \ge 6$ there exists $c=c_s$ such that $r_3(TP_s, n) > c n^2.$  The construction above has many copies of $TP_5$ so this leaves open the case $s=5$. Using the trivial lower bound $r_3(TP_4, n)$ we thus have $c_1 n^2/\log n < r_3(TP_5, n) < c_2 n^2.$

\begin{problem} [\cite{CMsparse}]
Determine the order of magnitude of $r_3(TP_5, n)$.
\end{problem}

The corresponding problems for $k$-graphs when $k>3$ are wide open.

\subsection{Cycles versus cliques}

For fixed $s \ge 3$ the graph  Ramsey number $r(C_s, n)=r(C_s, K_n)$ has been extensively studied. The case $s=3$ is one of the oldest questions in Ramsey theory and it is known that
$r(C_3,K_n)=\Theta(n^2/\log n)$ (see~\cite{AKS, Kim} and \cite{BK2, FGM} for recent improvements).
 The next case $r(C_4, K_n)$ seems substantially more difficult.
An old open problem of Erd\H os~\cite{E84} asks whether there is a positive $\epsilon$ for which
$r(C_4, K_n) = O(n^{2-\epsilon})$.  The current best upper bound $r(C_4,K_n) = O(n^2/\log^2 n)$ is an unpublished result of Szemer\'edi which was reproved in~\cite{ram-caro00} and the
current best lower bound is $\Omega(n^{3/2}/\log n)$ from~\cite{BK}. For longer cycles, the best known bounds can be found in~\cite{BK, SU}, and the order of magnitude of $r(C_s, K_n)$ is not known for any fixed $s \ge 4$.

There are several natural ways to define a cycle in hypergraphs. The two that have been investigated the most are tight cycles and loose cycles.

\subsubsection{Loose cycles versus cliques}

For $s \ge 3$, the loose cycle $LC_s^{(k)}$ is the $k$-graph with vertex set $\ZZ_{(k-1)s}$ and edge set $\{e_1, e_2, \ldots, e_s\}$ where $e_i=\{i(k-1)-k+2, \ldots, i(k-1)+1\}$. In other words, consecutive edges intersect in exactly one vertex and nonconsecutive edges are pairwise disjoint.
As usual, write $r_k(LC_s, n)$ for $r(LC_s^{(k)}, K_n^{(k)})$.  Since loose cycles are $k$-partite, it is easy to see that $r_k(LC_s, n)$ has polynomial growth rate for fixed $k,s$ so the question here is to determine the correct power of $n$.

\begin{theorem} [Kostochka-Mubayi-Verstra\"ete~\cite{kmvr}]\label{main}
There exists $c_1,c_2> 0$ such that for all $n \geq 1$,
$$\frac{c_1 n^{3/2}}{(\log n)^{3/4}} \leq r_3(LC_3,n) \leq c_2 n^{3/2}.$$
For $k \ge 3$, we also have $r_k(LC_3,n)=n^{3/2+o(1)}$.
\end{theorem}
Analogous to the basic result $r(3, n)= O(n^2/\log n)$ due to Ajtai, Koml\'{o}s and Szemer\'{e}di~\cite{AKS}, the authors conjectured something similar for hypergraphs.

\begin{conjecture}[\cite{kmvr}]\label{conjl6}
For all fixed $k \ge 3$, we have $r_k(LC_3,n) = o(n^{3/2})$.
\end{conjecture}

 Define the 3-graph $F=\{abc, abd, cde\}$. Cooper and Mubayi~\cite{CM} proved the following weaker version of Conjecture~\ref{conjl6} in the case $k=3$:
\begin{equation} \label{cm}r_3(\{LC_3, F, K_4^{(3)}-e\}, n) = O\left(
\frac{n^{3/2}}
{(\log n)^{1/2}}\right).\end{equation}
Notice that the three forbidden 3-graphs in (\ref{cm}) are all types of triangles, comprising three edges that cyclically share a vertex. Conjecture \ref{conjl6}
asks that we forbid only one of these three triangles, the loose triangle.

For longer cycles, the following general lower bounds were proved in~\cite{kmvr} which improve the bounds given by the standard probabilistic deletion method:

$$r_k(LC_s, n) > n^{1+1/(3s-1)+o(1)}.$$
Furthermore, there exists $c=c_k$ such that
$$r_k(LC_5, n)> c\left(\frac{n}{\log n}\right)^{5/4}.$$

\begin{conjecture} [\cite{kmvr}]
 For each $k \ge 3$, there exists $c=c_k$ such that
$r_k(LC_5, n)< c \, n^{5/4}$.
\end{conjecture}

M\'eroueh~\cite{Me} recently proved that $r_3(LC_5, n ) < c n^{4/3}$ and more generally that
$$r_3(LC_s, n) < c_s n^{1+\frac{1}{\lfloor (s+1)/2 \rfloor}}.$$
He also proved that for odd $s\ge 5$ and $k \ge 4$, $r_k(LC_s, n) < c_{k,s} n^{1+1/\lfloor s/2 \rfloor}$ which slightly improved  the exponent $1+1/(\lfloor s/2 \rfloor-1)$ proved by Collier-Cartaino, Graber and Jiang~\cite{CGJ} for all $k \ge 3$ and $s \ge 4$.

\subsubsection{Tight cycles versus cliques}
 For $k\ge 2$ and $s > 3$, the {\it tight cycle} $TC_s^{(k)}$ is the $k$-graph
with vertex set $\ZZ_s$ (integers modulo $s$) and edge set
$$\{\{i, i+1,  \ldots, i+k-1\}: i \in \ZZ_s\}.$$  We can view the vertex set of $TC_s^{(k)}$ as $s$ points on a circle and the edge set as the $s$  subintervals each containing $k$ consecutive vertices.

  When $s \equiv 0$ (mod 3) the tight cycle $TC_s^{(3)}$ is 3-partite, and in this case it is trivial to observe that $r_3(TC_s, n):= r(TC_s^{(3)}, K_n^{(3)})$ grows polynomially in $n$. The growth rate of this polynomial is not known for any $s>3$. When $s \not\equiv 0$ (mod 3) the Ramsey number is exponential in $n$.

  \begin{theorem}[Mubayi-R\"odl~\cite{MR,Mtc}]
   Fix $s \ge 5$ and $s \not\equiv 0$ (mod 3).  There are positive constants $c_1$ and $c_2$ such that
$$
2^{c_1n}<r_3(TC_s, n)<2^{c_2n^2\log n}.
$$

If $s \not\equiv 0$ (mod 3) and $s\ge 16$ or $s \in \{8,11,14\}$, there is a positive constant $c_s$ such that
$$r_3(TC_s, n) < 2^{c_s n \log n}.$$

  \end{theorem}

Note that when $s=4$, the cycle $TC_4^{(3)}$ is $K_4^{(3)}$ and in this case the lower bound was proved much earlier by Erd\H os and Hajnal~\cite{EH72}, and in fact has been improved to $2^{c_1 n \log n}$ more recently by Conlon, Fox and Sudakov~\cite{CFS}.

\begin{problem} Prove similar bounds for
$s\in\{ 4,5,7,10,13\}$ and determine whether the $\log$ factor in the exponent is necessary.
\end{problem}

 The problem of determining $r_k(TC_s, n):=r_k(TC_s^{(k)}, K_t^{(k)})$ for fixed $s>k>3$ seems harder as $k$ grows.
 It was shown in~\cite{MR} that we have a lower bound
$$r_k(TC_s, n) > 2^{c_{k, s} n^{k-2}}$$
The best upper bound that is known (for fixed $s>k$ and all $n$) is the trivial one $r_k(s, n)$.
Consequently, we have
$$2^{c_{k,s} n^{k-2}} < r_k(TC_s, n) <r_k(s, n)< \twr_{k-1}(n^{d_{s,k}}).$$
 Closing the  gap above seems to be a very interesting open problem. For the case $s=k+1$, one has a substantially better lower bound as
$$ r_k(TC_{k+1}, n) = r_k(k+1, n) > \twr_{k-2}(bn^{\log n})$$
where $b=b_k$.

\begin{problem} For fixed $s>k+1>4$ (in particular for $s=k+2$), determine whether
$r_k(TC_s, n)$ is   at least a tower function in a power of $n$ where the tower height grows with $k$.
\end{problem}

\section{Bounded degree hypergraphs}

Given a bounded degree graph $G$, the Ramsey number $r_2(G,G)$ has been studied extensively.  A famous result due to Chv\'atal, R\"odl, Szemer\'edi, and Trotter \cite{CRSTbd} says that if $G$ is a graph on $n$ vertices with maximum degree $\Delta$, then $r_2(G, G) \leq c_{\Delta}n$ where $c_{\Delta}$ depends only on $\Delta$.   This was later extended to 3-graphs by Cooly, Fountoulakis,
K\"uhn, Osthus \cite{CFK}, Nagle, Olsen, R\"odl, Schacht \cite{NOR}, and Ishigami \cite{Ish} independently, where the degree of a vertex $v$ in a hypergraph $H$ is the number of edges which contain $v$.  Using different methods, Cooly, Fountoulakis, K\"uhn, Osthus \cite{CFK2} and Conlon, Fox, Sudakov \cite{CFSbdk} extended this to general $k$-graphs.

\begin{theorem}[Cooly et al.~\cite{CFK2}, Conlon-Fox-Sudakov~\cite{CFSbdk}]

For all $\Delta, k \geq 1$, there is a $c(\Delta,k)$ such that for any $k$-graph $H$ on $n$ vertices with maximum degree $\Delta$, we have $$r_k(H,H) \leq c(\Delta, k) n.$$
\end{theorem}

In the special case that $H = LC^{(3)}_n$ or $TC_n^{(3)}$, the asymptotics for both $r_3(LC_n, LC_n)$ and $r_3(TC_n, TC_n)$
were determined in \cite{HLPRRSS} and \cite{HLPRRS} respectively.

\begin{theorem}[Haxell et al. \cite{HLPRRSS}]

$r_3(LC_n, LC_n) = \frac{5n}{2}(1 + o(1)).$

\end{theorem}

\begin{theorem}[Haxell et al. \cite{HLPRRS}]

For $n$ divisible by 3, we have $r_3(TC_n,TC_n) = \frac{4n}{3}(1 + o(1))$.  Otherwise, if $n$ is not divisible by 3, we have $r_3(TC_n,TC_n) = 2n(1 +o(1))$.

\end{theorem}

For $q$-colors, Gy\'arf\'as and Raeisi~\cite{GR} proved that
$$q+5 \le  r_3(LC_3; q) \le  3q+1.$$
It seems to be an interesting open problem to determine which bound above is closer to the truth. The lower bound appears more likely to be the answer.
There has been some work on determining $r_3(P; q)$ where $P=\{123, 345, 567\}$ is the 3-uniform loose path of length three.  In particular, the  lower bound $r_3(P; q)\ge q+6$ for all $q \ge 3$  is sharp for all $q \le 9$ (see~\cite{GR, PR}). The general upper bound which comes from the Tur\'an number of $P$ is again $3q+1$. This upper bound was improved by Luczak and Polcyn~\cite{LP} first to $(2+o(1))q$ and more recently to $\lambda q + O(\sqrt q)$ where the constant $\lambda=1.97466..$ is the solution to a particular cubic equation.

\section{Ordered Hypergraph Ramsey Problems}\label{ordsect}

In this section, we discuss several Ramsey-type results for ordered hypergraphs.  An \emph{ordered} $N$-vertex $k$-graph $H$ is a hypergraph whose vertex set is $[N] = \{1,\ldots, N\}$.  Given two ordered $k$-graphs $G$ and $H$ with vertex set $[n]$ and $[N]$ respectively, we say that $H$ \emph{contains} $G$ if there is a function $\phi:[n]\rightarrow [N]$ such that $\phi(i) < \phi(j)$ for all $1\leq i < j \leq n$, and $(v_1,\ldots, v_k) \in E(G)$ implies that $(\phi(v_1),\ldots, \phi(v_k)) \in E(H)$.  Given $q$ ordered $k$-graphs $H_1,\ldots, H_q$, the \emph{ordered Ramsey number} $\overline{r}_k(H_1,\ldots, H_n)$ is the minimum integer $N$, such that every $q$-coloring of the edges of the complete $k$-graph with vertex set $[N]$, contains a copy of $H_i$ in the $i$th color.

\subsection{Tight-paths and cliques in hypergraphs}

An \emph{ordered tight path} $P_s^{(k)}$ is an ordered $k$-graph with vertex set $[s]$, whose edges are of the form $(i , i + 1,\ldots, i + k - 1)$, for $1 \leq i \leq s - k + 1$.  The \emph{length} of an ordered tight path $P_s^{(k)}$ is the number of edges it contains, that is, $s-k + 1$.  In order to avoid the excessive use of superscripts, we write $P_s = P_s^{(k)}$ when the uniformity is already implied.   Two famous theorems of Erd\H os and Szekeres in \cite{ES35}, known as the monotone subsequence theorem and the cups-caps theorem, imply that $\overline{r}_2(P_s , P_n) = (n - 1)(s - 1) + 1$ and $\overline{r}_3(P_s,P_n) = {n+s - 4 \choose s-2} + 1.$   In \cite{FPSS}, Fox, Pach, Sudakov, and Suk extended their results to $k$-graphs and determined the correct tower growth rate for $\overline{r}_k(P_s,P_s)$.  Their results gave a geometric application related to the Happy Ending Theorem.\footnote{The main result in \cite{ES35}, known as the Happy Ending Theorem, states that for any positive integer $n$, any sufficiently large set of points in the plane in general position has a subset of $n$ members that form the vertices of a convex polygon.}  A few years later, Moshkovitz and Shapira~\cite{MS} sharpened the bounds for $\overline{r}_k(P_s,P_s)$ by determining an exact formula for $\overline{r}_3(P_{s},\ldots, P_{s})= \overline{r}_3(P_{s} ;q)$ with $q$ colors.

 \begin{theorem}[Moshkovitz-Shapira~\cite{MS}]
 Let $P_{q-1}(s)$ denote the number of $s\times \cdots \times s$ $(q-1)$-dimensional partitions with entries $\{0,1,\ldots, s\}$.  Then

 $$\overline{r}_3(P_s;q) = P_{q-1}(s) + 1.$$

 \end{theorem}

Soon after, Milans-Stolee-West~\cite{MSW} obtained an exact formula for $\overline{r}_k(P_{s_1}, \ldots, P_{s_q})$ for all $k, q \ge 2$, and $s_i \ge k$
 (see~\cite{CS,DLR} for some related results).

 \begin{theorem} [Milans-Stolee-West~\cite{MSW}]  Let $k, q \ge 2$, and $s_i > k$ for all $i \in [q]$. Let $J_1$ be the poset comprising disjoint chains $C_1, \ldots, C_q$, with $|C_i|=s_i-k$ for  $i \in [q]$
 and for $i\ge 1$, let $J_{i+1}$ be the poset whose elements are the ideals (down sets) of $J_i$ with order defined by containment.  Then $$\overline{r}_k(P_{s_1}, \ldots, P_{s_q})=|J_k|+1.$$
 \end{theorem}

Just as before, we will use the simpler notation $\overline{r}_k(P_s,n) = \overline{r}_k(P^{(k)}_s,K^{(k)}_n)$, and note that there is only one ordered complete hypergraph $K^{(k)}_n$ up to isomorphism.  Interestingly, the proof of the Erd\H os-Szekeres monotone subsequence theorem \cite{ES35} (see also Dilworth's Theorem \cite{D50}) actually implies that $\overline{r}_2(P_s , n) = (n - 1)(s - 1) + 1.$  For $k\geq 3$, estimating $\overline{r}_k(P_s, n)$ appears to be more difficult.  Clearly we have

\begin{equation} \label{easy} \overline{r}_k(P_s, n) \leq r_k(s,n)   \leq \twr_{k - 1}(O(n^{s - 2}\log n)).\end{equation}

In \cite{MS15}, the authors established the following connection between the ordered Ramsey number $\overline{r}_k(P_s,n)$ and the classical multi-color Ramsey number $r_k(n;q)$

\begin{theorem}[Mubayi-Suk \cite{MS15}]\label{mainx}
Let $k\geq 2$ and $s  \geq k + 1$.  Then for $q = s - k + 1$, we have $$r_{k -1}(\lfloor n/q\rfloor ;q)
\le \overline{r}_k(P_s, n) \le
r_{k -1}(n;q)
.$$

\end{theorem}

The upper bound in Theorem \ref{mainx} follows from the following argument.  Let $q = s - k +1$, $N = r_{k-1}(n;q)$, and suppose $\chi$ is a red/blue coloring on the $k$-tuples of $[N]$. We can assume $\chi$ does not produce a red tight-path of length $q$, since otherwise we would have a red $P_s$ and be done.  We define the coloring $\phi:{[N]\choose k-1} \rightarrow \{0,1,\ldots, q-1\}$ on the $(k-1)$-tuples of $[N]$, where $\phi(i_1,\ldots, i_{k-1}) = j$ if the longest red tight-path ending in vertices $(i_1,\ldots, i_{k-1})$ has length $j$.  Since $N = r_{k-1}(n;q)$, by Ramsey's theorem, we have a monochromatic clique of size $n$ in color $j$ for some $j \in \{0,1,\ldots, q-1\}$.  However, this clique would correspond to a blue clique with respect to $\chi$.  For the lower bound, set $N = r_{k -1}(\lfloor n/q\rfloor ;q) - 1$, and let $\chi$ be a $q$ coloring on the $(k-1)$-tuples of $[N]$ with colors $\{1,2,\ldots, q\}$, such that $\chi$ does not produce a monochromatic clique of size $\lfloor n/q\rfloor$.  Then let $\phi:{[N] \choose k}\rightarrow \{\textnormal{red,blue}\}$ such that for $i_1 < \cdots < i_k$, $\phi(i_1,\ldots, i_k)$ is red if and only if $\chi(i_1,\ldots, i_{k-1}) < \chi(i_2,\ldots, i_k)$.  It is easy to see that $\phi$ does not produce a red tight-path $P_s$.  With a slightly more complicated argument, one can show by contradiction that $\phi$ also does not produce a monochromatic blue clique of size $n$.

The arguments above can be easily extended to obtain the following result for multiple colors \cite{MS15}.  Let $k\geq 2$ and $s_1,\ldots, s_t \geq k + 1$.  Then for $q = (s_1 - k + 1)\cdots (s_t - k + 1)$, we have $$r_{k -1}(\lfloor n/q\rfloor ;q)
\le \overline{r}_k(P_{s_1}, \ldots, P_{s_t}, n) \le
r_{k -1}(n;q)
.$$

Together with known bounds for $r_{k-1}(n;q)$, Theorem \ref{mainx} has several consequences. First, we can considerably improve the upper bound for $\overline{r}_k(P_s, n)$ in (\ref{easy}) to $\overline{r}_k(P_s, n) \leq \twr_{k-1}(O(sn\log s)).$  In the other direction, the authors in \cite{MS15} showed that $\overline{r}_3(P_4,n) > 2^{cn}$, and Theorem \ref{mainx} implies that for $k \geq 4$ and $n>3k$,
\begin{enumerate}

\item $\overline{r}_k(P_{k + 3}, n) \geq \twr_{k-1}(cn),$

\item $\overline{r}_k(P_{k + 2}, n) \geq \twr_{k - 1}(c\log^2 n),$

\item  $\overline{r}_k(P_{k + 1}, n) \geq \twr_{k - 2}(cn^2).$

\end{enumerate}

We conjecture the following strengthening of the Erd\H os-Hajnal conjecture.

\begin{conjecture}\label{conj3}
For $k\geq 4$ fixed, $\overline{r}_k(P_{k + 1}, n)  \ge \twr_{k-1}(\Omega(n)).$
\end{conjecture}

For $s = k + 1$ in Theorem \ref{mainx}, we have $r_{k-1}(\lfloor n/2\rfloor, \lfloor n/2\rfloor) \leq \overline{r}_k(P_{k + 1},n) \leq r_{k-1}(n,n)$.  Hence, we obtain the following corollary which relates $\overline{r}_4(P_{5}, n)$ to the diagonal Ramsey number $r_3(n,n)$.

\begin{corollary}[\cite{MS15}]\label{equiv}
Conjecture \ref{3conj} holds if and only if
there is a constant $c>0$ such that
$$\overline{r}_4(P_5, n)  \ge 2^{2^{cn}}.$$
\end{corollary}

For the case when the size of $P_s$ tends to infinity and the size of $K_n$ is fixed, the first author in \cite{Mes} showed that $\overline{r}_3(P_s,4)< s^{21}$, and more generally for each $k \ge 3$, there  exists  $c>0$ such that for $s$ large, $$\twr_{k-2}(s^{c}) < \overline{r}_{k}(P_s,k+1) < \twr_{k-2}(s^{62}).$$

Unfortunately much less is known about $\overline{r}_k(P_s,k+2)$.  The main open problem here is to prove that $\overline{r}_3(P_s, 5)$ has polynomial growth rate, and more generally, that $\overline{r}_3(P_s, n)$ has polynomial growth rate for all fixed $n > 4$. The corresponding results for higher uniformity follow easily from the case $k=3$.

We next consider a version of the  Erd\H os-Hajnal hypergraph Ramsey problem with respect to tight-paths.

\begin{definition}   For integers $2\le k < s <n$ and $2 \le t \le {s \choose k}$, let $\overline{r}_k(s,t;P_n)$ be the minimum $N$ such that every red/blue coloring of the $k$-sets of $[N]$ results in a monochromatic blue copy of $P_n$ or has a set of $s$ vertices which induces at least $t$ red edges.
\end{definition}

Of course, $\overline{r}_k(s, {s \choose k}; P_n)= \overline{r}_k(s, P_n)$.  We will focus our attention on the smallest case $s=k+1$.  The following conjecture which parallels the Erd\H os-Hajnal conjecture for cliques was posed in~\cite{Mes}.

\begin{conjecture} [\cite{Mes}] \label{EHpaths}
For $3 \le t \le k$, there are positive $c =c(k,t)$ and $c'=c'(k,t)$ such that
$$\twr_{t-2}(n^c) < \overline{r}_k(k+1, t; P_n) < \twr_{t-2}(n^{c'}).$$
\end{conjecture}

This conjecture seems more difficult than the original problem of
Erd\H os and Hajnal. The current best lower bound is only an exponential function; unfortunately the constructions used for Theorem~\ref{MSEH} fail.
Standard arguments yield  an upper bound of the form $\twr_{t-1}(n^{c})$ for Conjecture~\ref{EHpaths}. This upper bound was improved in~\cite{Mes} to $\twr_{t-2}(n^{c})$.
Some further minor progress towards Conjecture~\ref{EHpaths} was made in~\cite{Mes} for the cases $t=3$ and $(k,t)=(4,4)$.

\subsection{Ordered $\ell$-power paths in graphs}

As mentioned above, The proof of Dilworth's theorem shows that $\overline{r}_2(P_s,n) = \overline{r}(P_s, P_n)=(s-1)(n-1)+1$.  On the other hand, we know that the classical Ramsey number $r_2(n,n)$ grows exponentially in $\Theta(n)$. One can consider
the case of ordered graphs that are denser than  paths but sparser than cliques.

\begin{definition} Given $\ell \ge 1$, the $\ell$th power $P_s^{\ell}$ of a path $P_s$ has ordered vertex set $v_1<\cdots  < v_s$ and edge set $\{v_iv_j: |i-j|\le \ell\}$. In particular, $P_s^1=P_s$.  The ordered Ramsey number $r(P_s^{\ell}, P_n^{\ell})$ is the minimum $N$ such that every red/blue coloring of ${[N] \choose 2}$ results in a red copy of $P_s^{\ell}$ or a blue copy of $P_n^{\ell}$.
\end{definition}

In~\cite{Mes} it was shown that the problem of determining $\overline{r}(P_n^{\ell}, P_n^{\ell})$ is closely related to the hypergraph ordered Ramsey function $\overline{r}_3(s, P_n)$.
 Conlon-Fox-Lee-Sudakov~\cite{CFLSord} asked whether $\overline{r}(P_n^{\ell}, P_n^{\ell})$ is polynomial in $n$ for every fixed $\ell\ge 1$.  Actually, the problem in~\cite{CFLSord} is about the Ramsey number of ordered graphs with bandwidth at most $\ell$ but $P^{\ell}_n$ contains all such graphs so an upper bound for $P^{\ell}_n$ provides an upper bound for the bandwidth problem. This question was answered by
Balko-Cibulka-Kr\'al-Kyn\v cl~\cite{BCKK}. Later a better bound was proved in \cite{Mes} for $\ell=2$.

\begin{theorem}[Balko-Cibulka-Kr\'al-Kyn\v cl~\cite{BCKK} ($\ell\ge 3$), Mubayi~\cite{Mes} ($\ell=2$)] \label{main}
There is an absolute constant $c>0$ and for every $\ell>0$ there exists $c=c_{\ell}$ such that
\begin{align}
\overline{r}(P_n^{\ell}, P_n^{\ell}) <
\begin{cases}
c\, n^{19.487} & \hbox{for } \ell=2 \\
c_{\ell} \, n^{128 \ell} & \hbox{for } \ell\ge 3.
\end{cases}
\end{align}
\end{theorem}

The main open problem here is to improve the exponents above. To our knowledge, there are no nontrivial lower bounds published for this problem.

\begin{problem} [Balko-Cibulka-Kr\'al-Kyn\v cl~\cite{BCKK}]
Determine the growth rate of $\overline{r}(P_n^{\ell}, P_n^{\ell})$ for every fixed $\ell \ge 2$.
\end{problem}

\section{A bipartite hypergraph Ramsey problem of Erd\H os}

We end with an old problem of Erd\H os that was perhaps posed to gain a better understanding of the growth rate of the diagonal Ramsey numbers.

\begin{definition} Let $S_{a,b} = (U,V,E)$ be the 3-graph with vertex set $U\cup V$, where $|U|=a$ and $|V| = b$, such that $E(S_{a,b}) = \{(x,y,z): x \in U \textnormal{ and } y,z \in V\}$. Write $S_n:=S_{n,n}$.
\end{definition}

 An old result due to Erd\H os (see \cite{E90}) says that $r_3(S_n,S_n) = 2^{O(n^2)}$, which is tight up to a constant factor in the exponent by the standard probabilistic method. We were not able to find a published proof of this result and we therefore present a proof below (of a stronger result).

 \begin{theorem}
For every $c>0$ and sufficiently large $n$,
$$r_3(S_n, S_n) < r_3(S_{2^{cn}, n}, S_{2^{cn}, n})<2^{3 n^2}.$$
 \end{theorem}
 \proof We begin with the simple observation that $r_3(S_{1,n}, S_{1,n}) < 1+r_2(n,n) < 4^n$. Indeed, if $r=1+r_2(n,n)$ and ${[r] \choose 3}$ is 2-colored by $\chi$ then we have an induced 2-coloring $\chi'$ of ${[r-1] \choose 2}$ where $\chi'(ij)=\chi(ijr)$. Because $r-1 = r_2(n,n)$ we have a monochromatic $n$-set under $\chi'$ and this yields a monochromatic $S_{1,n}$ under $\chi$ with $U=\{r\}$.

  Now we use a simple supersaturation trick to prove the result. Suppose that $N=2^{c' n^2}$ and $\chi$ is a 2-coloring of ${[N] \choose 3}$. For every $r$-set of $[N]$, where $r=4^n$, there is a monochromatic copy of $S_{1,n}$ in $\chi$. Hence the number of monochromatic copies of $S_{1,n}$ in $\chi$ is at least
  $$\frac{{N \choose r}}{{N-n-1 \choose r-n-1}} = \frac{(N)_{n+1}}{(r)_{n+1}}.$$
At least half of these monochromatic copies of  $S_{1,n}$ have the same color, say blue. Now, to each of these blue copies of $S_{1,n}$ with parts $|U|=1$ and $|V|=n$, we associate the $n$-set $V$.
A short calculation and the fact that $n$ is large shows that
    $$\frac{(N)_{n+1}}{(r)_{n+1}} > \left( \frac{ 2^c N  e  }{n} \right)^n > 2^{cn} {N \choose n}.$$
    Consequently, by the pigeonhole principle, there are at least $2^{cn}$ blue copies of $S_{1,n}$ associated to the same $n$-set $V$. These blue copies together form a blue copy of $S_n$ as desired.
 \qed

  Erd\H os stated that an important and difficult problem is to decide if his result can be strengthened to imply all triples that meet both $U$ and $V$.

  \begin{definition} Let $B_n = (U,V,E)$ be the 3-graph with vertex set $U\cup V$, where $|U| = |V| = n$, such that $E(B_n) = \{(x,y,z): x,y \in U, z\in V \textnormal{ or } x,y \in V, z\in U\}$.
\end{definition}

   Clearly we have
   $$2^{cn^2} < r_3(B_n,B_n) \leq r_3(n,n) \leq 2^{2^{c' n}}$$
where the lower bound follows from the probabilistic method.

\begin{problem}[Erd\H os \cite{E90}]  Improve the upper or lower bounds for $r_3(B_n, B_n)$.
\end{problem}

\end{document}